\documentclass[12pt]{amsart}

\begin{document}
\title[DRAZIN SPECTRA]{ DRAZIN SPECTRA OF BANACH SPACE OPERATORS AND BANACH ALGEBRA ELEMENTS}
\author[ENRICO BOASSO ]{ENRICO BOASSO}

\begin{abstract} Given a Banach Algebra $A$ and $a\in A$, several 
relations among the Drazin spectrum of $a$ and the Drazin spectra
of the multiplication operators $L_a$ and $R_a$ will be stated.
The Banach space operator case will be also examined.
Furthermore, a characterization of the Drazin spectrum
will be considered.\par
\end{abstract}
\maketitle

\noindent \bf{1. Introduction}\rm \vskip.3cm

\indent The main objective of the present article is to study the relationship 
between the notions of Drazin invertibility of Banach algebra elements and of Banach space operators. In particular, in the third section, after having recalled
several preliminary facts in section 2, given a unital Banach algebra $A$ and $a\in A$,  it will be proved that the Drazin spectrum of $a$
coincides with the Drazin spectrum of both $L_a$ and $R_a$, the left and right multiplication operators
defined by $a$ respectively. 
Furthermore, it will be proved that the Drazin spectrum of $a\in A$ coincides both with the union of the
right Drazin spectra of $L_a$ and $R_a$, and with the union of the descent spectra of the same
operators.  \par

\indent In addition, in the fourth section the relations among the ascent, the descent, the left and the right Drazin spectra of
a Banach space operator $T$ and the corresponding spectra of the multiplication
operators $L_T$ and $R_T$ will be also examined.\par 

\indent Finally, in section 5 a characterization of the Drazin spectrum will be considered. 
Actually, several well-known results concerning the Drazin spectrum of
Banach space operators will be extended to the frame of Banach algebras.\par
\vskip.3cm
\noindent \bf{2. Preliminary Definitions and Results}\rm 
\vskip.3cm
\indent From now on, $X$ will denote a Banach space and $L(X)$ will be the
Banach algebra of all operators defined on and with values in $X$. 
In addition, if $T\in L(X)$,  then $N(T)$ and $R(T)$ will stand for the null space and the
range of $T$ respectively. Moreover, $X^*$ will denote the dual space of $X$.
Recall that the \it descent \rm and the \it ascent \rm of $T\in L(X)$ are
$d(T) =\hbox{ inf}\{ n\ge 0\colon  R(T^n)=R(T^{n+1})\}$ and
$ a(T)=\hbox{ inf}\{ n\ge 0\colon  N(T^n)=N(T^{n+1})\}$
respectively, where if some of the above sets is empty, its infimum is then defined as $\infty$, see  
for example [19,6,9] and [5].\par
   
 \indent On the other hand, $A$ will denote a unital Banach algebra and its unit element will be  denoted by $e$. 
If $a\in A$, then $L_a \colon A\to A$
and $R_a\colon A\to A$ will denote the maps defined by
left and right multiplication respectively, namely $L_a(x)=ax$ and $R_a(x)=xa$,   
where $x\in A$. Concerning the null spaces and ranges of these maps, the following notation will be used:
$N(L_a)= a^{-1}(0)$, $R(L_a)=aA$, $N(R_a)= a_{-1}(0)$, and $R(R_a)=Aa$.\par

\indent Next follow the key notions of the present work.
 Given a Banach algebra $A$, an element $a\in A$ will be called
\it Drazin invertible\rm, if there exists a necessarily unique $b\in A$ and some $m\in \Bbb N$ such that
$$
a^mba=a^m, \hskip.5truecm  bab=b,\hskip.5truecm ab=ba.
$$
\indent If the Drazin inverse of $a$ exists, then it will be denoted by $a^D$. In addition,
the \it index \rm of $a$, which will be denoted by  $index (a)$, is the least non-negative integer $m$ for which the
above equations hold. When $index(a)=1$,
$a$ will be said \it group invertible\rm, and in this case its Drazin inverse will be refered as the group inverse
of $a$, which will be denoted by $a^{ \sharp}$,
see [12,10,15,20] and [6]. \par

\indent On the other hand, the notion of \it regularity \rm was introduced in [16]. In order to learn the definition
and the main properties of the aforesaid notion, see [16] and [19]. Recall that
given a unital Banach algebra $A$ and a regularity $\mathcal R\subseteq A$,
the \it spectrum derived from the 
regularity \rm $ \mathcal R$  is defined as $\sigma_{\mathcal R} (a)=\{ \lambda\in \Bbb C\colon a-\lambda\notin  \mathcal R\}$,
where $a\in A$ and $ a-\lambda$ stands for  $a-\lambda e$, see [16]. In addition, the \it resolvent set of $a$ defined by the regularity $\mathcal R$ \rm
is the set $\rho_{\mathcal R} (a)=\{ \lambda\in \Bbb C\colon a-\lambda\in  \mathcal R\}$. Next consider the set
$\mathcal{DR} (A)$ = $\{ a\in A\colon  \hbox{  a is Drazin invertible}\}$. According to
Theorem 2.3 of [6], $\mathcal{DR} (A)$ is a regularity. This fact led to the following definition,
see [6].\par

\newtheorem*{def1}{Definition 1} 
\begin{def1} Let $A$ be a unital Banach algebra. The Drazin spectrum
of an element $a\in A$ is the set
$$
\sigma_{ \mathcal{ DR} }(a) = \{\lambda\in \Bbb C \colon a-\lambda\notin  \mathcal{ DR}(A)\}.
$$
\end{def1}

\indent Naturally $\sigma_{ \mathcal{ DR} }(a)\subseteq \sigma(a)$, and according to Theorem 1.4 of [16], the Drazin spectrum
of a Banach algebra element satisfies the spectral mapping theorem for analytic functions
defined on a neighbourhood of the ususal spectrum which are non-constant on each component of its 
domain of definitioin, see also Corollary 2.4 of [6]. In addition, according to Proposition 2.5 of [6],
$\sigma_{ \mathcal{ DR} }(a)$ is a closed subset of $\Bbb C$. \par

\indent When $A=L(X)$, $X$ a Banach space, the left and the right Drazin spectra
of an operator have been introduced, see Definition 2.5 of [7] and Definition 1.1 of [3,4],
see also [1,2]. Before recalling these notions, consider
the sets 
\begin{align*}
\mathcal{ LD}(X)&=\{ T\in L(X)\colon a(T) \hbox{ is finite and  }
R(T^{ a(T) +1 })\hbox{ is closed}\}, \\
\mathcal{ RD}(X)& = \{ T\in L(X)\colon d(T) \hbox{  is finite and  }
R(T^{ d(T) })\hbox{ is closed}\}.\\
\end{align*}
\indent Recall that, according to [19, pp.134 and 138], $\mathcal{ LD}(X)$ and 
$\mathcal{ RD}(X)$ are regularities. \par

\newtheorem*{def2}{Definition 2} 
\begin{def2} Let $X$ be a Banach space. An operator $T\in L(X)$ will be called left
Drazin invertible (respectively right Drazin invertible), if $T\in \mathcal{ LD}(X)$ (respectively if
 $T\in \mathcal{ RD}(X)$). Given $T\in L(X)$, the left Drazin spectrum
of $T$ (respectively the right Drazin spectrum of $T$) is defined as the spectrum derived from the
regularity  $\mathcal{ LD}(X)$ (respectively from the regularity $\mathcal{ RD}(X)$). 
These spectra will be denoted by $\sigma_{  \mathcal{ LD} }(T)$ and $\sigma_{  \mathcal{ RD} }(T)$ 
respectively, and $\rho_{\mathcal LD}(T)$ and $\rho_{\mathcal RD}(T)$ will stand for
the corresponding resolvent sets.
\end{def2}

\indent In the conditions of Definition 2, according again to Theorem 1.4
of [16], $\sigma_{  \mathcal{ LD} }(T)$ and $\sigma_{  \mathcal{ RD} }(T)$ 
satisfy the spectral mapping theorem under the same hypothesis
that the Drazin spectrum does. 
 Moreover, according to [19, pp.137 and 139],  $\sigma_{  \mathcal{ RD} }(T)$ and $\sigma_{  \mathcal{ LD} }(T)$
are closed subsets of $\Bbb C$. See also the recent article [5], where  the properties of the spectrum  $\sigma_{  \mathcal{ LD} }(T)$
were intensively studied.\par

\indent In order to introduce the remaining spectra that will be used in this work,
the classes of Banach space operators of finite descent and of finte ascent will be considered.\par

\indent Let $X$ be a Banach space. Recall that, according to [19], the sets 
\begin{align*}
\mathcal R^a_4(X) &= \{ T\in L(X)\colon d(T)\hbox{  is finite} \},\\
 \mathcal R^a_9(X) &= \{ T\in L(X)\colon a(T) \hbox{  is finite}\},\\
\end{align*}
are regularities with well-defined  spectra derived from them, namely the \it descent spectrum \rm and the \it ascent spectrum \rm 
respectively. The descent spectrum was also studied in [9] and it will be denoted by
 $\sigma_{desc}(T)$, $T\in L(X)$. On the other hand, $\sigma_{asc}(T)$ will stand for the ascent spectrum.
Note that $\sigma_{desc}(T)\subseteq\sigma_{ \mathcal{ RD} }(T) $ and $\sigma_{asc}(T)\subseteq\sigma_{ \mathcal{ LD} }(T) $.
\par

\indent The following theorem is a consequence of well-known facts. However, since in this theorem
 the relationships among all the spectra above recalled 
will be presented, the complete proof is given. \par

\newtheorem*{theo3}{Theorem 3} 
\begin{theo3} Let $X$ be a Banach space and consider $T\in L(X)$. Then
$$
 \mathcal{ DR}(L(X)) = \mathcal{ LD}(X)\cap \mathcal{ RD}(X)
=  \mathcal R^a_9(X)\cap \mathcal R^a_4(X).
$$
In particular,
$$
\sigma_{ \mathcal{ DR} }(T) =  \sigma_{ \mathcal{ LD} }(T) \cup  \sigma_{ \mathcal{ RD} }(T)
 =  \sigma_{asc }(T) \cup \sigma_{desc}(T).
$$
\begin{proof} Since all the spectra involved in the statement
correspond to regularities, according to [19, p.130], 
it is enough to prove the first statement. \par

\indent Note that, according to Theorem 4 of [15], 
$$\mathcal R^a_4(X)\cap \mathcal R^a_9(X) =\mathcal{ DR}(L(X)).$$

\indent Furthermore, according again to Theorem 4 of [15], it is clear that $\mathcal{ LD }(X) \cap \mathcal{ RD }(X) $
$\subseteq\mathcal{ DR}(L(X))$.\par

\indent  On the other hand, if $T$ has a Drazin inverse,
then according to Theorems 3 and 4 of [15], the ascent and descent of
$T$ are finite, $k$ = $index (T)$ = $a(T)$ = $d(T)$, and $X=N(T^k)\oplus R(T^k)$.
Since for $n\ge d(T)$, $R(T^n)=R(T^{ d(T) })$ is closed,  $T\in\mathcal{ LR}(X)\cap  \mathcal{ RD}(X)$. 
Therefore, $\mathcal{ DR}(L(X))$ = $\mathcal{ LD}(X)\cap\mathcal{ RD}(X)$.
\end{proof}
\end{theo3}
\vskip.3cm
\noindent \bf{3. The Drazin Spectrum in Banach Algebras}\rm 
\vskip.3cm
\indent In order to prove the main results of this section, first of all
a characterization of Drazin invertible Banach algebra elements will be presented. On the other hand, given $A$ a Banach algebra and $a\in A$, recall that, according to Proposition 4 of [8, section 5, Chapter 1], 
$\sigma(a)=\sigma(L_a)$. Note, however, that the identity $\sigma(a)=\sigma(R_a)$ can be also obtained using the same 
argument of the mentioned reference. 
In the following theorem a similar relationship will be also proved for the Drazin spectrum.\par 

\newtheorem*{theo4}{Theorem 4} 
\begin{theo4} 
\indent Let $A$ be a unital Banach algebra. Then, the following statements are equivalent.\par
\begin{align*}
&(i)& &\hbox{The element }a\in A \hbox{ is Drazin invertible and }index (a)=k,&\\
&(ii)& &L_a\in L(A) \hbox{ is Drazin invertible and }index (L_a)=k,\\
&(iii)& &R_a\in L(A) \hbox{ is Drazin invertible and }index (R_a)=k.\\
\end{align*}
\indent Moreover, in this case, if $b$ is the Drazin inverse of $a$,
then $L_b$ (respectively $R_b$) is the Drazin inerse of $L_a$ (respectively $R_a$).\par

\indent As a consequence,\par
\indent \hskip.3truecm (iv) \hskip1truecm $\sigma_{\mathcal{DR} }(a) =\sigma_{\mathcal{DR}}(L_a)=\sigma_{\mathcal{DR}}(R_a)$.\par

\begin{proof} In first place, suppose that $a\in A$ is group invertible. It is then clear that $L_a$ and $R_a$ have group inverses. In addition,
if  $a^{\sharp}$ is the group inverse of $a$, then $(L_a)^{\sharp} =L_{a^{\sharp}}$ and $(R_a)^{\sharp} = R_{a^{\sharp}}$.  \par
\indent On the other hand, if $W\in L(A)$ is the group inverse of
$L_a$, then 
$$
L_a=(L_a)^2W= L_{a^2}W.
$$
In particular, $aR=a^2R$. What is more, since 
$$
L_a= W(L_a)^2= WL_{a^2},
$$
$a^{-1}(0)=(a^2)^{-1}(0)$. Therefore, according to Proposition 7 of [14], 
$a$ is group invertible, and according to what has been proved, $W=(L_a)^{\sharp} =L_{a^{\sharp}}$.\par

\indent Next consider the general case. As before, it is clear that
if $a^D$ is the Drazin inverse of $a\in A$, then $L_{a^D}$ is the Drazin inverse of $L_a\in L(A)$.
Moreover, $m=index(L_a)\leq index(a)=k$.\par
\indent On the other hand, since $(L_a)^m=L_{a^m}$, and
$$
(L_a)^m =(L_a)^m L_{a^D} L_a,\hskip.5truecm L_{a^D}L_aL_{a^D}=L_{a^D},\hskip.5truecm L_{a^D}L_a=L_aL_{a^D},  
$$
a straightforward calculation proves that
$$
a^m=a^ma^Da,\hskip.5truecm a^Daa^D=a^D,\hskip.5cm a^Da=aa^D.
$$
\indent In particular, $k\leq m$.\par
\indent Suppose that $L_a\in L(A)$ is Drazin invertible and let $k=index(L_a)$.  
According to Corollary 5 of [20], $(L_a)^k=L_{a^k}$ is group invertible in $L(A)$.
However, according to what has been proved, 
$a^k$ is group invertible in $A$, which, according again to Corollary 5 of [20],
implies that $a$ is Drazin invertible and $index (a)=k$. \par
\indent The equivalence between (i) and (iii) can be proved in a similar way.\par
\indent The last statement is a consequence of the statements (i) -(iii).
\end{proof}
\end{theo4}

\indent Note that, according to Theorem 4 of [15], Corollary 2.13 of  [5]
is equivalent to the identity $\sigma_{\mathcal{DR} }(T)=\sigma_{\mathcal{DR} }(L_T)$, for $T\in L(X)$, 
$X$ a Banach space.\par
\indent In addition, note  that the equivalence of the statements (i)-(iii) of Theorem 4 also holds for group
invertible elements of rings with units, see the proof of Theorem 4 and Proposition 7 of [14].\par 

\indent The next theorem will show that the Drazin spectrum can be computed
using the descent and the right Drazin spectra of the multiplication operators.\par

\newtheorem*{theo5}{Theorem 5} 
\begin{theo5} Consider a unital Banach algebra $A$, and let $a\in A$. Then,
the following statements hold.\par
\begin{align*} 
&(i)&  &\sigma_{asc}(R_a)\subseteq \sigma_{desc}(L_a), \sigma_{asc}(L_a)\subseteq \sigma_{desc}(R_a),& \\
&(ii)& &\sigma_{\mathcal{DR}}(a)= \sigma_{\mathcal{RD}}(L_a)\cup \sigma_{\mathcal{RD}}(R_a) =\sigma_{desc}(L_a)\cup \sigma_{desc}(R_a).\\
\end{align*} 
\begin{proof} 
\indent The first statement is a consequence of Proposition 3.2 (a) of [21]. \par
\indent Concerning the second statement, according to Theorem 3 and  Theorem  4 (iv), it is clear that
$$
 \sigma_{desc}(L_a)\cup \sigma_{desc}(R_a)\subseteq \sigma_{\mathcal{RD}}(L_a)\cup \sigma_{\mathcal{RD}}(R_a) \subseteq
\sigma_{\mathcal{DR}}(a).
$$

\indent On the other hand, suppose that $\lambda\notin (\sigma_{desc}(L_a)\cup \sigma_{desc}(R_a))$. Then,
according to Proposition 3.2 (c) of [21], Theorem 4 of [15], and  Theorem 4 (iv), $\lambda\notin \sigma_{\mathcal{DR}}(a)$.

\end{proof}
\end{theo5}

\vskip.3cm
\noindent \bf{4. Drazin Spectra of Banach Space Operators}\rm \vskip.3cm

\indent In a Banach algebra $A$ the notion of Drazin invertibility can be considered both for
elements of the algebra and for operators defined on the Banach space $A$. This two 
different applications of the notion under consideration have been compared in the previous section. 
Moreover, Theorem 5 relates the Drazin spectrum with the descent and the
right Drazin spectra of two operators. On the other hand,
the concepts of ascent, descent, and left and right Drazin invertibility  
apply only to bounded linear maps. Therefore, in order to relate them
to their algebraic versions,  a particular class of algebras should be considered,
namely $A=L(X)$, where $X$ is a Banach space. In what follows these relationships
will be studied. However, in first place some preparation is needed.\par

\newtheorem*{rem6}{Remark 6} 
\begin{rem6}\rm Let $A$ be a Banach algebra. Recall that an element
$a\in A$ is called \it regular\rm, if there is $b\in A$ such that
$a=aba$. For example, if $a$ is Drazin invertible and $k=index(a)$,
then $a^k$ is regular. When $A=L(X)$, $X$ a Banach space,
it is well known that necessary and sufficiente for $T\in L(X)$ to be regular
is the fact that $N(T)$ and $R(T)$ have direct complements in $X$,
see Theorem 1 of [10].
In particular, if $H$ denotes a Hilbert space, $T\in L(H)$ 
is a regular operator if and only if $R(T)$ is closed.\par

\indent Recall in addition that in a $\Bbb C^*$-algebra,
necessary and sufficient for an element $a\in A$ to be regular
is the fact that the range of $L_a\in L(A)$ is closed,
naturally an equivalent formulation can be proved using
$R_a\in L(A)$ instead of $L_a$, see Theorems 2 and 8 of [13].
Observe in particular that $a\in A$ is regular if and
only if $L_a$ (respectively $R_a$) is a regular
operator on $L(A)$.
\par   

\indent On the other hand, it is well-known that if 
$A$ and $B$ are operators defined on a Hilbert
space $H$, necessary and sufficient for $R(A)\subseteq R(B)$
is the fact that there exists $C\in L(H)$ such that
$A=BC$, see Theorem 1 of [11]. What is more,
the same result can be obtained for operators defined on 
a Banach space $X$,
provided that $N(B)$ has a direct complement in $X$,
see the proof of Theorem 1 of [11] and [11, p. 415]. In what follows, a similar formulation
regarding the null spaces instead of the ranges will 
be considered.\par

\indent Let $X$ be a Banach space, and let $A$ and $B$
belong to $L(X)$. Suppose that $B$ is regular and
$N(A)=N(B)$. Then, there exists $C\in L(X)$ such that 
$A=CB$.\par
\indent In fact, according to Theorem 1 of  [10], there are $M$ and $N$ two closed subspaces of $X$
such that 
$$
N(A)\oplus M=N(B)\oplus M  =X,\hskip.5truecm R(B)\oplus N=X.
$$
\indent Now well, since $\tilde{A}=A\mid_M^{R(A)}\colon M\to R(A)$
and $\tilde{B}=B\mid_M^{R(B)}\colon M\to R(B)$ are two isomorphic
Banach space operators, it is possible to consider the Banach
space operator $\tilde{C}=\tilde{A}\tilde{B}^{-1}\colon R(B)\to R(A)$.
Next define $C\in L(X)$ as follows: $C\mid_{R(B)}= \iota_{R(A), X} \tilde{C}\colon R(B)\to X$,
and $C\mid_N\colon N\to X$ any bounded and linear map, where
$\iota_{R(A), X} \colon R(A)\to X$ denotes the inclusion map.
Then, it is clear that $A=CB$.\par
\end{rem6}

\indent In the following proposition some basic facts will be presented.
Two of them were proved before, however, since they will be central
for the main results of this section, they will be listed and the
original references will be indicated. Recall that if $X$ is a Banach space,
then $L_T$ and $R_T\in L(L(X))$ denote the left and right multiplication
operators defined by $T\in L(X)$ respectively.\par 
\par

\newtheorem*{prop7}{Proposition 7} 
\begin{prop7}Consider a Banach space $X$, and let $T\in L(X)$.
\par
\begin{align*} 
&(i)& &\hbox{If } desc(L_T) \hbox{ is finite, then } desc(T) \hbox{ is finite. In addition, }desc(T)\\
&\hbox{  }& &\le desc (L_T).&\\
&(ii)& &\hbox{If } desc(T)=d \hbox{ is finite and } N(T^{d+1}) \hbox{ has a direct complement,}\\
&\hbox{  }& &\hbox{then }
desc(L_T) \hbox{ is finite. Moreover, } desc(L_T)=desc (T).&\\
\end{align*}
\begin{align*}
&(iii)& &\hbox{Necessary and suffcient for } asc (L_T) \hbox{ to be finite, is the fact that }\\
&\hbox{  }& & asc(T) \hbox{ is finite. Furthermore, in this case } asc (L_T)= asc(T).&\\
&(iv)& &\hbox{If } R(L_T) \hbox{ is closed, then } R(T) \hbox{  is closed. } \\
&(v)& &\hbox{If } R(T) \hbox{ is closed and  }N(T) \hbox{ has a direct complement, then}\\
&\hbox{  }& &R(L_T)\hbox{ is closed.}&\\
&(vi)& &\hbox{If } desc(R_T)  \hbox{ is finite, then }asc(T) \hbox{ is finite. In addition, }asc(T)\\
&\hbox{  }& &\le desc (R_T).&\\
&(vii)& &\hbox{If } asc(T)=a \hbox{ is finite and } T^{a+1} \hbox{ is a regular operator,} \hbox{ then }&\\
&\hbox{  }& &desc(R_T)\hbox{ is finite. Moreover, } desc(R_T)=asc (T).&\\
&(viii)& &\hbox{If } asc(R_T)=a  \hbox{ is finite and there exists $k\ge a$ such that } &\\
&\hbox{  }& & R(T^{k+1}) \hbox{ has a direct complement, then } desc (T) \hbox{ is finite. }&\\
&\hbox{  }& &\hbox{Furthermore, }desc(T)= asc(R_T).&\\
&(ix)& &\hbox{If } desc(T)   \hbox{ is finite, then  }  asc(R_T)  \hbox{ is finite. What is more, }&\\
&\hbox{  }& &asc(R_T)\le desc(T).&\\
\end{align*}
\begin{proof} The first statement was proved in [9, p. 263].\par
\indent Suppose that $desc(T)$ is finite. Then there is a natural number $d$
such that $R(T^d)=R(T^{d+1})$. If in addition $N(T^{d+1})$ has a direct
complement, then according to the proof of Theorem 1 of [11], see also
[11, p. 415], there is $C\in L(X)$ such that $T^d=T^{d+1}C$,
which implies that the descent of $L_T$ is finite and $desc(L_T)\le desc(T)$.
According to the first statement, $desc(L_T)= desc(T)$.\par

\indent  In order to prove the third statement, suppose that there is a positive integer $a$ such that
$N((L_T)^{a+1})=N((L_T)^a)$, equivalently $N(L_{T^{a+1}})=N(L_{T^a})$.
Consider $x\in N(T^{a+1})$ and construct $P=P^2\in L(X)$ such that
$R(P)=<x>$, the vector space generated by the element $x$.
Since $x\in N(T^{a+1})$ and $X=<x>\oplus N(P)$, it is clear that
$P\in N(L_{T^{a+1}})=N(L_{T^a})$. In particular, $L_{T^a}(P)=0$.
However, $0=T^aP(x)=T^a(x)$, which implies that $x\in N(T^a)$.
What is more, $asc (T)\le asc (L_T)$.\par

\indent On the other hand, suppose that there is a natural number $a$ such that
$N(T^{a+1})=N(T^a)$. Let $S\in N(L_{T^{a+1}})$. Since necessary and sufficient for $S$ to belong to $N(L_{T^{a+1}})$
is that $R(S)\subseteq N(T^{a+1})=N(T^a)$, $S\in N(L_{T^a})$, equivalently $N(L_{T^{a+1}})=N(L_{T^a})$.
As a consequence,  $asc(L_T)\le asc (T)$, and according to what
has been proved, $asc(T)=asc(L_T)$.\par  

\indent Next suppose that $R(L_T)$ is closed. Consider a Cauchy sequence $(T(x_n))_{n\in\Bbb N}\subseteq X$,
and let $z\in X$ and $f\in X^*$ such that $f(z)=\parallel f\parallel =1$. For each $n\in \Bbb N$
define the operator $U_n\in L(X)$ as follows: $U_n(x)= T(x_n)f(x)$. It is not difficult to prove that $(U_n)_{n\in\Bbb N}\subseteq L(X)$
is a Cauchy sequence such that $(U_n)_{n\in\Bbb N}\subseteq R(L_T)$.
Therefore, there exists $S\in L(X)$ such that $TS=\lim_{n\to\infty} U_n$.
In particular, $T(S(z))=\lim_{n\to\infty}T(x_n)$.\par

\indent In order to prove the fifth statement, consider a sequence of operators
$(S_n)_{n\in \Bbb N}\subseteq L(X)$ and an operator $Q\in L(X)$
such that $(TS_n)_{n\in \Bbb N}$ converges to $Q$ in $L(X)$. In particular,
if $x\in X$, then $(TS_n(x))_{n\in \Bbb N}$ converges to $Q(x)$.
According to the fact that $R(T)$ is closed, it is clear that
$R(Q)\subseteq R(T)$. Now well, since $N(T)$ has a direct complement,
according again to the proof of Theorem 1 of [11], see also
[11, p. 415], there is an operator $U\in L(X)$ such that 
$Q=TU$, in particular, $Q\in R(L_T)$ and $R(L_T)$ is closed.\par 

\indent On the other hand, suppose that $desc(R_T)$ is finite. Note that
necessary and suffcient for $desc(T)$ to be finite is the fact that
there is $S\in L(X)$ such that $T^d=ST^{d+1}$,
which clearly implies that $N(T^{d+1})\subseteq N(T^d)$. In particular,
$asc(T)$ is finite. Furthermore, $asc(T)\le desc(R_T)$.\par

\indent If $N(T^{a+1})= N(T^a)$ and $T^{a+1}$ is a regular operator,
then, according to Remark 6, there is $S\in L(X)$ such that
$T^a=ST^{a+1}$. Consequently, $desc(T)$ is finite. Moreover, 
according to what has been proved, $desc(T)=asc(T)$.\par

\indent Next suppose that $asc (R_T)=a$. Then,
$N(R_{T^a})=N(R_{T^n})$, for all $n\ge a$. Consider 
$P=P^2\in L(X)$ such that $R(P)=R(T^{k+1})$. 
Since $k\ge a$ and $(I-P)(R(T^{k+1}))=0$, $I-P\in N(R_{T^{k+1}})
=N(R_{T^a})$. Consequently, $(I-P)(R(T^a)) =0$,
which implies that $R(T^a)\subseteq N(I-P)=R(T^{k+1})\subseteq
R(T^a)$. Since $k\ge a$, $desc(T)=a$.\par

\indent In order to prove the last statement, suppose there is
a natural number $d$ such that $desc(T)=d$, equivalently
$R(T^d)=R(T^{d+1})$. Consider $S\in L(X)$ such that 
$S\in N(R_{T^{d+1}})$. Since $0=S(R(T^{d+1}))=S(R(T^d))$,
$S\in N(R_{T^{d}})$. In particular $asc (R_T)\le  d=desc(T)$. 
\end{proof}

\end{prop7}

\indent Note that the equivalence between the finitness of $asc(T)$ and
$asc(L_T)$ considered in Proposition 7 (iii) was proved in Proposition 2.12 of [5]. 
The reason for giving a different proof of this results lies in the fact that the identity
$asc(T)=asc(L_T)$ can be also obtained.\par

\newtheorem*{theo8}{Theorem 8} 
\begin{theo8} Consider a Banach space $X$, and $T\in L(X)$. Then, the 
following statements hold.\par
\begin{align*}
&(i)&&\sigma_{desc}(T)\subseteq \sigma_{desc}(L_T), \hbox{ }\sigma_{asc}(T)=\sigma_{asc}(L_T),&\\
&(ii)&&\sigma_{  \mathcal{RD}}(T)\subseteq  \sigma_{ \mathcal{RD}}(L_T), \hbox{      }\hbox{ }\sigma_{\mathcal{LD}}(T)\subseteq  \sigma_{\mathcal{LD}}(L_T),&\\
\end{align*}
\begin{align*}
&(iii)&&  \sigma_{asc}(T)\subseteq  \sigma_{desc}(R_T), \hbox{ } \sigma_{asc}(R_T)\subseteq  \sigma_{desc}(T).&\\
\end{align*} 
\indent Furthermore,
\begin{align*}
&(iv)& &\sigma_{desc}(L_T)\setminus\sigma_{desc}(T)\subseteq \sigma_{asc}(T), \hbox{ }\sigma_{desc}(L_T)\cap\rho_{asc}(T)=\sigma_{desc}(T)&\\
&\hbox{ }&&\cap\rho_{asc}(T),&\\
&(v)& &\sigma_{ \mathcal{RD}}(L_T)\setminus\sigma_{ \mathcal{RD}}(T)\subseteq \sigma_{ \mathcal{LD}}(T), \hbox{ }\sigma_{ \mathcal{RD}}(L_T)\cap\rho_{ \mathcal{LD}}(T)=\sigma_{ \mathcal{RD}}(T)&\\
&\hbox{ }&&\cap\rho_{ \mathcal{LD}}(T),&\\
&(vi)& &\sigma_{ \mathcal{LD}}(L_T)\setminus\sigma_{ \mathcal{LD}}(T)\subseteq \sigma_{ \mathcal{RD}}(T), \hbox{ }\sigma_{ \mathcal{LD}}(L_T)\cap\rho_{ \mathcal{RD}}(T)=\sigma_{ \mathcal{LD}}(T)&\\
&\hbox{ }&&\cap\rho_{ \mathcal{RD}}(T),&\\
&(vii)& &\sigma_{desc}(R_T)\setminus\sigma_{asc}(T)\subseteq \sigma_{des}(T), \hbox{ }\sigma_{desc}(R_T)\cap\rho_{desc}(T)=\sigma_{asc}(T)&\\
&\hbox{ }&&\cap\rho_{desc}(T),&\\
&(viii)& &\sigma_{desc}(T)\setminus\sigma_{asc}(R_T)\subseteq \sigma_{desc}(R_T), \hbox{ }\sigma_{desc}(T)\cap\rho_{desc}(R_T)=\sigma_{asc}(R_T)&\\
&\hbox{ }&&\cap\rho_{desc}(R_T).&\\
\end{align*}

\begin{proof}  As regard the first statement, it is a consequence of statements (i) and (iii) of Proposition 7.\par
\indent In order to prove the second statement, apply statement (iv) of Proposition 7 to what has been proved. \par
\indent In addition, the third statement is a consequence of statements (vi) and (ix) of 
Proposition 7.\par
\indent According to Theorem 3 and Theorem 5 (ii), $\sigma_{desc}(L_T)\setminus\sigma_{desc}(T)\subseteq \sigma_{asc}(T)$.
On the other hand, according to the first statement, it is clear that  
$$
\sigma_{desc}(T)\cap\rho_{asc}(T)\subseteq\sigma_{desc}(L_T)\cap\rho_{asc}(T).
$$
The other inclusion is a consequence of the following fact.\par
\begin{align*}
\sigma_{desc}(L_T)\cap\rho_{asc}(T) \setminus \sigma_{desc}(T)\cap\rho_{asc}(T)&=(\sigma_{desc}(L_T) \setminus \sigma_{desc}(T))\\
\cap\rho_{asc}(T)&\subseteq \sigma_{asc}(T)\cap\rho_{asc}(T)=\emptyset.\\
\end{align*}
\indent Statements (v) - (viii) can be proved in a similar way.
\end{proof}
\end{theo8}

\indent In the following theorem the left Drazin spectrum of left and right multiplication operators 
will be used to compute the Drazin spectrum of Hilbert bounded and linear maps. Compare this 
result with Theorem 5 (ii).\par

\newtheorem*{theo9}{Theorem 9} 
\begin{theo9} Consider a Hilbert space $H$, and let $T\in L(H)$. Then, the following statements hold.\par
\begin{align*} 
&(i)&&\sigma_{desc}(T)=\sigma_{desc}(L_T).&\\
&(ii)&&\sigma_{  \mathcal{RD}}(T)= \sigma_{ \mathcal{RD}}(L_T).&\\
&(iii)&&\sigma_{\mathcal{LD}}(T)= \sigma_{\mathcal{LD}}(L_T).&\\
&(iv)&&  \sigma_{\mathcal{LD}}(T)=  \sigma_{\mathcal{RD}}(R_T).&\\
&(v)&&\sigma_{\mathcal{RD}}(T)=  \sigma_{\mathcal{LD}}(R_T).&\\
\end{align*} 
\indent As a result,\par
\begin{align*} 
&(vi)& &\sigma_{\mathcal{DR}}(T)=\sigma_{\mathcal{LD}}(L_T)\cup\sigma_{\mathcal{LD}}(R_T).&
\end{align*} 
\begin{proof}  The first statement is a consequence of Theorem 8 (i) and Proposition 7 (ii).\par

\indent According to Theorem 8 (ii), in order to prove the second statement, it is enough to show 
that  if $0\notin \sigma_{  \mathcal{RD}}(T)$,
then $0\notin \sigma_{ \mathcal{RD}}(L_T) $.
Consider $T\in L(H)$ such that  $d=desc(T)$ is finite. Then, 
according to Proposition 7 (ii), $d=desc(L_T)$. In addition,
since $R(T^d)$ is closed, according to Proposition 7 (v),
$R(L_{T^d})$ is closed. \par

\indent According again to Theorem 8 (ii),  
the proof of (iii) can be concluded, if the condition 
$0\notin \sigma_{  \mathcal{LD}}(T)$ implies that
$0\notin \sigma_{ \mathcal{LD}}(L_T) $.
Suppose then that  $T\in L(H)$ is such that  $a=asc(T)$ is finite. Then, according to Proposition 7 (iii), $asc(L_T)=a$. In addition,
according to Proposition 7 (v), since $R(T^{a+1})$ is closed, $R(L_{T^{a+1}})$ is also
closed. \par

\indent Next consider $T\in L(H)$ such that $0\notin \sigma_{\mathcal{RD}}(R_T)$.
Then, according to Proposition 7 (vi),
$a=asc(T)$ is finite and $a\le desc(R_T)$. 
Furthermore,
if $d=desc (R_T)$, since, for $n\ge d$, $R(R_{T^n})=R(R_{T^d})$ is closed,
according to Theorem 8 of [13] and  Theorem 
1 of [10], $R(T^n)$ is closed for all $n\ge d$. However, according to
Lemma 7 of [19], $R(T^{a+1})$ is closed. In particular, $0\notin \sigma_{\mathcal{LD}}(T)$.\par 

\indent Now consider $T\in L(H)$ such that $0\notin \sigma_{\mathcal{LD}}(T)$.
Let $a=asc(T)$. Since $R(T^{a+1})$ is closed, according 
to Theorem 1 of [10], $T^{a+1}$ is a regular operator. Then,
according to Proposition 7 (vii), $d=desc(R_T)$ is finite, actually $d=desc(R_T)=asc(T)=a$.
Furthermore, according to Theorem 2 of [13], $R(R_{T^d})=R(R_{T^{a+1}})$ is closed.
Therefore, $0\notin \sigma_{\mathcal{RD}}(R_T)$.\par 

\indent In order to prove the fifth statement, consider $T\in L(H)$ such that $0\notin  \sigma_{\mathcal{LD}}(R_T)$,
and let $a=asc(R_T)$. Since $R(R_{T^{a+1}})$ is closed, according to Theorem 8 of [13],
$T^{a+1}$ is a regular operator, which, according to Theorem 1 of [10],
 implies that $R(T^{a+1})$ is closed. Consequently, according to Proposition 7 (viii), $desc(T)$ is finite,
actually $desc(T)=asc(R_T)=a$. However, since  $R(T^a)=R(T^{a+1})$, $0\notin \sigma_{\mathcal{RD}}(T)$.\par

\indent On the other hand, if $0\notin \sigma_{\mathcal{RD}}(T)$, then $d=desc(T)$ is finite,
and according to Proposition 7 (ix), $asc (R_T)$ is finite, what
is more, $a=asc(R_T)\le desc(T)=d$.
In addition, since $R(T^{d+1})$ is closed, according to Theorem 1 of [10] and  Theorem 2 of
[13], $R(R_{T^{d+1})})$ is closed. Therefore, according to Lemma 7 of [19],
$R(R_{T^{a+1}})$ is closd, equivalently $0\notin  \sigma_{\mathcal{LD}}(R_T)$.\par

\indent The last statement is a consequence of Theorem 5 (ii), and what has been proved.
\end{proof}
\end{theo9}

\indent In the same conditions of Theorem 9, note that inclusion $ \sigma_{desc}(L_T)\subseteq \sigma_{desc}(T)$ was proved in [9, p. 265]
using Theorem 1 of [11]. In fact, Proposition 7 (ii) consists in a generalization
of this argument to Banach spaces. In addition, Theorem 9 (iii) was proved in Proposition 2.12
of [5].\par
\vskip.3cm

\noindent \bf{5. A Characterization of the Drazin Spectrum}\rm \vskip.3cm

\indent In what follows  a description of the Drazin spectrum of a Banach algebra element will be presented. However, first of all
some notation is needed. Let $K\subseteq \Bbb C$ a compact set. Then $iso (K)$ will
stand for the set of all isolated points of $K$ and $acc (K)=K\setminus iso (K)$. \par

\newtheorem*{rem10}{Remark 10} 
\begin{rem10}\rm
\indent Consider a unital Banach algebra $A$ and $a\in A$. As in the case of a  
Banach space operator, the resolvent function of 
$a$, $R(\cdot, a)\colon\rho(a)\to A$, is holomorphic and $iso(\sigma(a))$ coincides with the set of
isolated singularities of $R(\cdot ,  a)$. Furthermore, as in the case of an operator,
see [22, p. 305], if $\lambda_0\in iso(\sigma(a))$, then it is possible to
consider the Laurent expansion of  $R(\cdot ,  a)$ in terms of $(\lambda-\lambda_0)$.
In fact,
$$
R(\lambda ,  a)=\sum_{n\ge 0}a_n(\lambda-\lambda_0)^n +\sum_{n\ge 1}b_n(\lambda-\lambda_0)^{-n},
$$
where $a_n$ and $b_n$ belong to $A$ and are obtained in an standard way using the
functional calculus. In addition, this representation is valid when $0<\mid \lambda-\lambda_0\mid<\delta$, for 
any $\delta$ such that $\sigma(a)\setminus \{\lambda_0\}$ lies on or outside the circle $\mid\lambda-\lambda_0\mid=\delta$. 
What is more important, the discussion of [22, pp. 305 and 306]
can be repeated for elements in a unital Banach algebra. Consequently, $\lambda_0$
will be called \it a pole of order $p$ of   $R(\cdot ,  a)$, \rm if there is $p\ge 1$ such that
$b_p\neq 0$ and $b_m=0$, for all $m\ge p+1$. The set of poles of $a$
will be denoted by $\Pi(a)$. \par

\indent On the other hand, define $\mathcal{IES}(a)=iso(\sigma(a))\setminus \Pi(a)$,
the set of all \it isolated essential singularities \rm of $R(\cdot, a)$. In the following
theorem $\Pi(a)$ and $\mathcal{IES}(a)$ will be compared with the corresponding sets of 
the operators $L_a$ and $R_a$. However, before going on two preliminary results
need to be considered.\par

\indent First, note the following equivalence. Let $U\subseteq \Bbb C$ be an open set, and consider an analytic function
$g\colon U\to A$. Then, it is not difficult to prove that $L_g\colon U \to L(A)$ is analytic and $(L_g)'(z)=L_{g'}(z)$.
Conversely, if $g\colon U\to A$ is a function such that $L_g\colon U \to L(A)$ 
is analytic, then, a direct calculation shows that $g\colon U \to A$ is analytic and $g'(z)=(L_g)'(z)(e)$, 
where $e$  is the identity of $A$.\par 

\indent Second, if $X$ is a Banach space operator and $T\in L(X)$, 
then the following identity holds.\par
\indent $(i)$\hskip3.1cm $\Pi (T)=\partial\sigma(T)\cap\rho_{\mathcal{DR}}(T)$.\par
\indent In fact, according to Theorem  5.8-A of [22], 
$$
\Pi (T)\subseteq \partial\sigma(T)\cap\rho_{\mathcal{DR}}(T).
$$

\indent To prove the other inclusion, according to Theorem 1.5 of [9], it is clear that
$$
\partial\sigma(T)\cap\rho_{\mathcal{DR}}(T)\subseteq\partial\sigma(T)\cap\rho_{desc}(T)=\Pi (T).
$$
\end{rem10}

\newtheorem*{theo11}{Theorem 11} 
\begin{theo11}  Let $A$ be a unital Banach algebra and consider $a\in A$. Then, the 
following statements hold.\par
\begin{align*}
&(i)& &\Pi(a)=\Pi (L_a)=\Pi(R_a),&\\
&(ii)& &\mathcal{IES}(a)=\mathcal{IES}(L_a)=\mathcal{IES}(R_a).&\\
\end{align*}
\indent In particular, if $X$ is a Banach space and $T\in L(X)$, then 
$$
\Pi(T)=\Pi(L_T)=\Pi(R_T), \hskip.3cm\mathcal{IES}(T)=\mathcal{IES}(L_T)=\mathcal{IES}(R_T).
$$
\begin{proof} According to Proposition 4 of [8, section 5, Chapter 1],  it is enough to prove the first statement. Moreover, according again to
Proposition 4 of [8, section 5, Chapter 1] and to Theorem 4 (iv), and Remark 10 (i), $\Pi (L_a)=\Pi (R_a)$.
Therefore, in order to conclude the proof, the identity $\Pi (a)=\Pi(L_a)$ will be proved.\par

\indent It is clear that
$$
R(\cdot, a)\colon \rho(a)\to A,\hbox{ and } R(\cdot, L_a)\colon \rho(a)\to L(A),
$$
are holomorphic functions.\par
\indent In addition, a straightforward calculation proves that  
$$
R(\cdot, L_a) =L_{R(\cdot, a)}\colon \rho(a)\to L(A).
$$

\indent Next consider $\lambda_0\in iso(\sigma (a))=iso (\sigma(L_a))$ a pole of order $p$ of
$R(\cdot, a)$.  Note that according Remark 10, $\lambda_0$ is a pole of order
$p\ge 1$ of $R(\cdot, a)$ if and only if $g(\lambda)=(\lambda-\lambda_0)^pR(\cdot, a)$ is an analytic function
in a neighborhood $U$ of $\lambda_0$ such that $g(\lambda_0)\neq 0$. Then, according to Remark 10, 
$L_g\colon U\to L(A)$
is a holomorphic function, and since $(L_g)(\lambda_0)(e)=g(\lambda_0)$, $(L_g)(\lambda_0)\neq 0$.
However, for $\lambda\in U\setminus \{\lambda_0\}$,
$$
L_g(\lambda)=(\lambda-\lambda_0)^pR(\lambda, L_a). 
$$
\indent Therefore, according to the equivalence of the previous paragraph applied to $R(\cdot, L_a)$, $L_a\in L(A)$,
$\lambda$ is a pole of order $p$ of $L_a$. In particular, $\Pi(a)\subseteq \Pi (L_a)$.\par
\indent Conversely, suppose that $\lambda_0$ is a pole of order $p$ of $L_a$. 
Consequently, as above, there is a $U$ a neighborhood of $\lambda_0$ such that $h\colon U\to L(A)$,
$h(\lambda)=(\lambda-\lambda_0)^pR(\lambda, L_a)$ is a holomorphic funtion such that
$h(\lambda_0)\neq 0$. However, since $R(\lambda, L_a) =L_{R(\lambda, a)}$,
$h(\lambda)=L_{(\lambda-\lambda_0)^pR(\lambda, a)}$, for $\lambda\in U$ and
$\lambda\neq \lambda_0$. Define $g\colon U\setminus \{\lambda_0\}\to A$,
$g(\lambda)=(\lambda-\lambda_0)^pR(\lambda, a)$. A straightforward
calculation shows that if $(\lambda_n)_{n\in\Bbb N}\subseteq U\setminus \{\lambda_0\}$
converges to $\lambda_0$, then $(g(\lambda_n))_{n\in\Bbb N}$ converges to $h(\lambda_0)(e)$,
where $e$ is the unit of $A$.
Extend the function $g\colon U\to A$ defining $g(\lambda_0)=h(\lambda_0)(e)$.
Now well, an easy calculation proves that  $h(\lambda_0)(b)=g(\lambda_0)b=(L_g)(\lambda_0)(b)$,
for any $b\in A$. Therefore, since for every $\lambda\in U$, $h(\lambda)=(L_g)(\lambda)$,
according to Remark 10, $g\colon U\to A$ is analytic. Moreover, since $h(\lambda_0)\neq 0$, $g(\lambda_0)\neq 0$. 
As a result, $\lambda_0$ is a pole of order $p$ of $R(\cdot ,a)$, that is $\Pi(L_a)\subseteq \Pi(a)$.
\end{proof}
\end{theo11}

\indent Thanks to Theorem 11  several well-known results concerning the Drazin spectrum of a Banach space
operator will be extended to Banach algebra elements, see Theorems 5.8-A and 5.8-D of [22]
and Theorem 12 of [10]. 
\par

\newtheorem*{theo12}{Theorem 12} 
\begin{theo12} Let $A$ be a unital Banach algebra and consider $a\in A$. Then,
the following statements hold.\par
\begin{align*}
&(i)&&\Pi(a)=\sigma(a)\cap\rho_{\mathcal{DR}}(a),  \sigma(a)=\sigma_{\mathcal{DR}}(a)\cup \Pi(a),&\\
&(ii) & &\sigma_{\mathcal{DR}}(a)\cap \Pi(a)=\emptyset, iso(\sigma(a))\cap \sigma_{\mathcal{DR}}(a)=\mathcal{IES}(a).&\\
&(iii)& &acc(\sigma(a))= acc(\sigma_{\mathcal{DR}}(a)), iso(\sigma_{\mathcal{DR}}(a))=\mathcal{IES}(a).&\\
&(iv)& &\sigma_{\mathcal{DR}}(a) =acc(\sigma (a))\cup \mathcal{IES}(a).&\\
&(v)& &\hbox{Necessary and sufficient for } a \hbox{ to be Drazin invertible is that }\\
&\hbox{ }&&\lambda=0 \hbox{ is a pole of the resolvent operator of } a.&\\
\end{align*}
\begin{proof} It is well-known that the first and the second statements hold for the particular case $A=L(X)$,
where $X$ is a Banach space, see for example Theorems 5.8-A and 5.8-D of [22],
Theorem 2.1 of [17], and Corollary 2.10 of [17] . Consequently, in order to prove (i) and (ii),
apply Proposition 4 of [8, section 5, Chapter 1], Theorem 4 (iv), Theorem 11, and what has been recalled.\par
\indent The third and the fourth statements are a consequence of what has been proved.
In order to prove the last statement, use Theorem 11 (i) and Theorem 12 of [10]. 
\end{proof}
\end{theo12}
\indent Note that Theorem 12 (iv) was proved in Proposition 1.5 of [18].\par
\newtheorem*{rem13}{Remark 13} 
\begin{rem13}\rm Let $X$ be a Banach space and consider $T\in L(X)$. Note that according
to the results of Banach space operators recalled in the proof of Theorem 12, a straightforward calculation
proves that the sets $\mathcal{IES}(T)$, $iso(\sigma(T))\cap\sigma_{\mathcal{LD}}(T)$, $iso(\sigma(T))\cap \sigma_{\mathcal{RD}}(T)$,
and $iso(\sigma(T))\cap \sigma_{desc}(T)$ coincide. When $A$ is a unital Banach algebra, similar presentations can be obtained for
$\mathcal{IES}(a)$ using the operators $L_a$ and $R_a\in L(A)$ and the descent and the left and the right Drazin
spectra of these operators.\par
\end{rem13}

\noindent \bf{Acknowledgements.} \rm The author wish to thank
Professors Pietro Aiena and Maria Burgos, for these researchers have kindly
sent to him several recent publications authored by them which are related to the subject of the present article.
The author would also express his indebtedness to the referees of this manuscript, for their observations
have improved the original version of the work. In particular the author is thankful to the
referee who suggested a shorter proof of Proposition 7 (iv) and pointed out to him [5].\par

\vskip.3truecm
\noindent Enrico Boasso\par
\noindent E-mail address: enrico\_odisseo@yahoo.it

\end{document}